\documentclass[11pt,twoside,leqno]{aomamlt2e} 
  \pageno{1335}
\received{November 9, 2003}
 
 \renewcommand{\thesection}{\arabic{section}}
\makeatletter
\renewcommand{\@seccntformat}[1]{\csname
the#1\endcsname.\hspace{0.5em}\setcounter{Subsec}{0}\setcounter{Subsubsec}{0}}\makeatother

\newcommand{\diml}[1]{{{#1}}}
\newcommand{\dima}[1]{{{#1}}}               % This removes all the color

\newcommand{\zima}[1]{{{#1}}}

\def\whp{w.h.p.}
\def\wupp{w.u.p.p.}

\def\ie{i.e.,\ }

\def\reals{\mathbb{R}}

\def\3gl{{\sc 3-gl}}

\newcommand{\one}{\mathbf{1}}

\newcommand{\bea}{\begin{eqnarray*}}
\newcommand{\eea}{\end{eqnarray*}}

\newcommand{\mat}{\left(\!\!\begin{array}{cc}}
\newcommand{\rix}{\end{array}\!\!\right)}

\newcommand{\h}{c}

\newcommand{\E}{\mathbb{E}}

\newcommand{\R}{\mathbb{R}}

\newcommand{\St}{{\cal S}}
\newcommand{\D}{\mathcal{D}}

\newtheorem{theorem}{Theorem}
\newtheorem{lemma}[theorem]{Lemma}

\newtheorem{proposition}[theorem]{Proposition}
\theorembodyfont{\upshape}

\newtheorem{definition}[theorem]{{\it Definition}}

\def\ex{{\mathbb E}}

\def\e{\varepsilon}

\begin{document}
\currannalsline{162}{2005} 

 \title{The two possible values of the\\ chromatic number of a random graph}
\acknowledgements{Work performed while the first author was at Microsoft
Research.}
 
\twoauthors{Dimitris Achlioptas}{Assaf Naor}

 \institution{Department of Computer Science, University of California Santa Cruz\\
\email{optas@cs.ucsc.edu}\\
\vskip-9pt
Microsoft Research, Redmond, WA
\\
\email{anaor@microsoft.com}}

 \shorttitle{The chromatic number of a random graph}

\centerline{\bf Abstract}
\vglue6pt 

\dima{Given} $d \in (0,\infty)$ let $k_d$ be the smallest integer
$k$ such that $d<2k\log k$. We prove that the chromatic number of
a random graph $G(n,d/n)$ is either $k_d$ or $k_d+1$
almost \dima{surely}.  

\vglue-24pt
\phantom{problem}
\section{Introduction}
\vglue-6pt

The classical model of random graphs, in which each possible edge on $n$
vertices is chosen independently with probability $p$, is denoted
by $G(n,p)$. This model, introduced by Erd\H{o}s and R\'enyi in
1960, has been studied intensively in the past four decades. We
refer to the books~\cite{alon-spencer}, \cite{bollobas}, \cite{jlr} and the
references therein for   accounts of many remarkable results on
random graphs, as well as  for their  connections to various
areas of mathematics. In the present paper we \dima{consider
random graphs of bounded average degree, \ie $p=d/n$ for some
fixed $d \in (0,\infty)$.}

One of the most important invariants of a graph $G$ is its
chromatic number $\chi(G)$, namely the minimum number of colors
required to color its vertices so that no pair of adjacent
vertices has the same color. Since the mid-1970s, work on
$\chi\left(G(n,p)\right)$ has been in the forefront of random
graph theory, motivating some of the field's most significant
developments. \dima{Indeed, one of the most fascinating facts
known~\cite{luczak2} about random \diml{graphs is} that for every
$d\in (0,\infty)$ there exists an integer $k_d$ such that almost
surely $\chi(G(n,d/n))$ is either $k_d$ or $k_d+1$.} \dima{The
value of $k_d$ itself, nevertheless,   remained a mystery}.

To date, the best known~\cite{luczak} estimate for
$\chi(G(n,d/n))$  confines it to an\break
\vglue-12pt\noindent interval of length about
$d\cdot \frac{29\log\log d}{2(\log d)^2}$. In our main result we
reduce this length to $2$. \zima{Specifically, we prove}
\begin{theorem}\label{thm:main}
Given $d \in (0,\infty)${\rm ,} let $k_d$ be the smallest integer $k$
such that $d<2k\log k$. With probability that tends to $1$ as
$n\to \infty${\rm ,}
\begin{eqnarray*}%\label{eq:luczak2}
\chi\left(G(n,d/n)\right) \in \{k_d, k_d+1\}
    \enspace .
\end{eqnarray*}
\end{theorem}

Indeed, we determine $\chi\left(G(n,d/n)\right)$ {\em exactly\/}
for roughly half of \dima{all $d \in (0,\infty)$.}
\begin{theorem}\label{thm:yuval}
If $d \zima{\in [(2}k-1) \log k,  2k \log k)${\rm ,} then with
probability that tends to $1$ as $n\to \infty${\rm ,}
$$
\chi\left(G(n,d/n)\right) = k+1.
$$
\end{theorem}

\diml{The} first questions regarding the chromatic number \diml{of
$G(n,d/n)$}  were raised in the original Erd\H{o}s-R\'enyi
paper~\cite{er} from 1960.  It was only until the 1990's,
though, that  any progress was made on the problem.  Specifically, by the mid
1970s, the expected value of \dima{$\chi(G(n,p))$} was known up to
a factor of two for the case of fixed $p$, \diml{due to the work
of Bollob{\'a}s and Erd\H{o}s~\cite{berd} and Grimmett and
McDiarmid~\cite{grimett}.} This gap remained in place for
\diml{another} decade until, in a celebrated paper,
Bollob{\'a}s~\cite{bela} \diml{proved that for every constant
$p\in (0,1)$, almost surely $\chi(G(n,p)) ={\frac{n}{2 \log
n}\log\left(\frac{1}{1-p}\right)(1+o(1))}$.} \L
uczak~\cite{luczak} later extended this result to \diml{all $p>
d_0/n$,} where $d_0$ is a universal constant. 

\diml{Questions regarding the concentration of the chromatic
number were first examined} in a seminal paper of Shamir and
Spencer~\cite{shamir} in the mid-80s. \diml{They} showed that
$\chi\left(G(n,p)\right)$ is concentrated in an interval of length
$O(\sqrt{n})$ for all $p$ and on an interval of length 5 for
$p<n^{-1/6-\e}$. \L uczak~\cite{luczak2} showed that, for
$p<n^{-1/6-\e}$ the chromatic number is, in fact, concentrated on
an interval of length 2. Finally, Alon and Krivelevich~\cite{alon}
extended $2$-value concentration to all $p<n^{-1/2-\e}$.

The Shamir-Spencer theorem mentioned above was based on analyzing
the so-called {\em vertex exposure martingale}. Indeed, this was
the first use of martingale methods in random graph theory. Later,
a much more refined martingale argument was the key step in
Bollob{\' a}s' evaluation of the asymptotic value of
$\chi(G(n,p))$. This influential line of reasoning has fuelled
many developments in probabilistic combinatorics --- in particular
all the results mentioned above~\cite{luczak}, \cite{luczak2}, \cite{alon} rely on
martingale techniques.

\dima{Our proof of Theorem~\ref{thm:main} is largely analytic,
breaking \zima{with more} traditional combinatorial arguments.}
The starting point for our approach is recent progress on the
theory of sharp thresholds. Specifically, using Fourier-analytic
arguments, Friedgut~\cite{Ehud} has obtained a deep criterion for
the existence\break of sharp thresholds for random graph properties.
Using Friedgut's theorem,\break Achlioptas and
Friedgut~\cite{frie_sharp} proved that the
probability that $G(n,d/n)$ is\break $k$-colorable drops from almost 1
to almost $0$ as $d$ \zima{crosses} an interval whose length tends
to $0$ with $n$. Thus, in order to prove that $G(n,d/n)$ is almost
surely $k$-colorable it suffices to prove that $\liminf_{n \to
\infty} \Pr[G(n,d'/n) \mbox{ is $k$-colorable}]\break>0$, for some
$d'>d$. To do that we use the second moment method, which is based
on the following special case of the Paley-Zygmund inequality: for
any nonnegative random variable $X$, $\Pr[X>0]\ge (\ex X)^2/\ex
X^2$.

Specifically, the number of $k$-colorings of a random graph is the
sum, over all $k$-partitions $\sigma$ of its vertices (into $k$
``color classes"), of the indicator that $\sigma$ is a valid
coloring. To estimate the second moment of the number of\break $k$-colorings we thus
need to understand the correlation between these indicators. It
turns out that this correlation is determined by $k^2$ parameters:
given two $k$-partitions $\sigma$ and~$\tau$, the probability that
both of them are valid colorings is determined by the number of
vertices that receive color $i$ in $\sigma$ and color $j$ in~$\tau$, where $1 \leq i,j \leq k$. 

In typical second moment arguments, the main task lies in using
probabilistic and combinatorial reasoning to construct a random
variable for which correlations can be controlled. \zima{We}
achieve this here \zima{by focusing on} the number, $Z$, of
$k$-colorings in which all color classes have exactly the same
size. However, we face an additional difficulty, of an entirely
different nature: the correlation parameter is inherently high
dimensional. %In particular,
\zima{As a result,} estimating $\ex Z^2$ reduces to a certain
entropy-energy inequality over $k \times k$ doubly stochastic
matrices and, thus, our argument shifts to the analysis of an
optimization problem over the Birkhoff polytope. Using geometric
and analytic ideas we establish the desired inequality as a particular
case of a general optimization principle that we formulate
(Theorem~\ref{thm:stochastic}). We believe that this principle
will find further applications, for example in probability and
statistical physics, as moment estimates are often characterized
by similar trade-offs.

\section{Preliminaries}\label{section:preliminaries}

We will say that a sequence of events ${\mathcal{E}}_n$ occurs
with high probability (w.h.p.) if $\lim_{n \to
\infty}\Pr[{\mathcal{E}}_n]=1$ and with uniformly positive
probability (w.u.p.p.) if $\liminf_{n \to
\infty}\Pr[{\mathcal{E}}_n]>0$. Throughout, we will consider $k$
to be arbitrarily large but fixed, while $n$ tends to infinity. In
particular, all asymptotic notation is with respect to $n \to
\infty$.

To prove Theorems~\ref{thm:main} and~\ref{thm:yuval} it will be
convenient to introduce a slightly different model of random
graphs. Let $G(n,m)$ denote a random (multi)graph on $n$ vertices
with precisely $m$ edges, each edge formed by joining two
vertices selected uniformly, independently, and with replacement.
\zima{The following elementary argument was first suggested by Luc
Devroye (see~\cite{Chv144}).}
\begin{lemma}\label{lem:luc} Define
$$
u_k\equiv \frac{\log k}{\log
k-\log(k-1)}<\left(k-\frac12\right)\log k \enspace .
$$
If $c>u_k${\rm ,} then a random graph $G(n,m=cn)$ is {\rm w.h.p.\/} \pagebreak
non\/{\rm -}\/$k$\/{\rm -}\/colorable.
\end{lemma}

\Proof Let $Y$ be the number of $k$-colorings of a random graph
$G(n,m)$. By Markov's inequality, $\Pr[Y>0] \leq \ex[Y] \leq k^n
\left(1-1/k\right)^{m}$ since, in any fixed $k$-partition a random
edge is monochromatic with probability at least $1/k$. \zima{For
$c > u_k$, we have $k(1-1/k)^c <1$, implying $\ex[Y] \to 0$.}
\Endproof\vskip4pt  
Define
$$
c_k \equiv k \log k \enspace .
$$
We will prove
\begin{proposition}\label{thm:threshold} If $c < \h_{k-1}${\rm ,} then a
random graph \zima{$G(kn,m=ckn)$} is {\rm \wupp}\ $k$\/{\rm -}\/colorable.
\end{proposition}

Finally, as mentioned in the introduction, we will use the
following result of~\cite{frie_sharp}.

\begin{theorem}[Achlioptas and Friedgut~\cite{frie_sharp}]\label{thm:sharp} 
\hskip-6pt Fix $d^*>d>0$. If $G(n,d^*/n)$ is $k$-colorable 
{\rm \wupp}\ then
$G(n,d/n)$ is $k$-colorable {\rm \whp}
\end{theorem}

\zima{We now prove Theorems~\ref{thm:main} and~\ref{thm:yuval}
given Proposition~\ref{thm:threshold}.}

\demo{Proof of Theorems~{\rm \ref{thm:main}}
and~{\rm \ref{thm:yuval}}} A random graph $G(n,m)$ may contain some
loops and multiple edges. Writing $q=q(G(n,m))$ for the number of
such blemishes we see that their removal results in a graph on $n$
vertices whose edge set is uniformly random among all edge sets of
size $m-q$. Moreover, note that if $m \leq cn$ for some constant
$c$, then \whp\ $q = o(n)$. Finally, note that the edge-set of a
random graph $G(n,p=2c/n)$ is uniformly random conditional on its
size, and that \whp\ this size is in the range $cn \pm n^{2/3}$.
Thus, if $A$ is any monotone decreasing property that holds with
probability at least $\theta>0$ in $G(n,m= cn)$, then $A$ must
hold with probability at least $\theta-o(1)$ in $G(n,d/n)$ for any
constant $d<2c$ and similarly, for increasing properties and $d>2c$.
Therefore, Lemma~\ref{lem:luc} implies that $G(n,d/n)$ is \whp\
non-$k$-colorable for $d\ge (2k-1) \log k>2u_k$.

To prove both theorems it thus suffices to prove that $G(n,d/n)$
is \whp\ $k$-colorable if $d<2\h_{k-1}$. Let $n'$ be the smallest
multiple of $k$ greater than $n$. Clearly, if $k$-colorability
holds with probability $\theta$ in $G(n',d/n')$ then it must hold
with probability at least $\theta$ in $G(t,d/n')$ for all $t \leq
n'$. Moreover, for $n \leq t\leq n'$, $d/n' = (1-o(1))d/t$ . Thus,
if $G(kn',m=ckn')$ is $k$-colorable \wupp, then $G(n,d/n)$ is
$k$-colorable \wupp\ for all $d<2c$. Invoking
Proposition~\ref{thm:threshold} and Theorem~\ref{thm:sharp} we
thus conclude that $G(n,d/n)$ is \whp\ $k$-colorable for all
$d<2c_{k-1}$.
\Endproof\vskip4pt  

\zima{In the next section we reduce the proof of
Proposition~\ref{thm:threshold} to an analytic inequality, which
we then prove in the remaining sections.}

\section{The second moment method and stochastic
matrices}\label{section:second}

In the following we will only consider random graphs $G(n,m=cn)$
where $n$ is a multiple of $k$ and $c>0$ is a constant. We will say that a partition of
$n$ vertices into $k$ parts is balanced if each part contains
precisely $n/k$ vertices. Let $Z$ be the number of balanced
$k$-colorings.
Observe that each balanced partition is a valid $k$-coloring with
probability $(1-1/k)^m$. Thus, by Stirling's approximation,
\begin{eqnarray}\label{eq:expectation}
\E Z=\frac{n!}{[(n/k)!]^k}\left(1-\frac{1}{k}\right)^{m}=
\Omega\left(\frac{1}{n^{(k-1)/2}}\right)
\left[k\left(1-\frac{1}{k}\right)^{c}\right]^n \enspace .
\end{eqnarray}
Observe that the probability that a $k$-partition is a valid
$k$-coloring is maximized when the partition is balanced.
Thus, focusing on balanced partitions reduces the number of colorings considered by only a polynomial
factor, while significantly simplifying calculations. We will show that $\E Z^2 < C \cdot (\E Z)^2$ for some
$C=C(k,c)<\infty$. By~\eqref{eq:expectation} this reduces to proving
$$
\E Z^2=
O\left(\frac{1}{n^{k-1}}\right)\left[k\left(1-\frac{1}{k}\right)^c\right]^{2n}
\enspace .
$$
This will conclude the proof of Proposition~\ref{thm:threshold}
since $\Pr[Z>0] \geq {(\ex Z)^2}/{\ex Z^2}$.
\smallbreak
Since $Z$ is the sum of ${n!}/{[(n/k)!]^k}$ indicator variables,
one for each balanced partition, we see that to calculate $\E Z^2$
it suffices to consider all pairs of balanced partitions and, for
each pair, bound the probability that both partitions are valid
colorings. For any fixed pair of partitions $\sigma$ and $\tau$,
since edges are chosen independently, this probability is the
$m$th power of the probability that a random edge is bichromatic
in both $\sigma$ and $\tau$. If $\ell_{ij}$ is the number of
vertices with color $i$ in $\sigma$ and color $j$ in $\tau$, this
single-edge probability is
\begin{eqnarray*}
  1-\frac{2}{k}+\sum_{i=1}^k\sum_{j=1}^k
  \left(\frac{\ell_{ij}}{n}\right)^2.
\end{eqnarray*}
Observe that the second term above is independent of the
$\ell_{ij}$ only because $\sigma$ and $\tau$ are balanced.

Denote by $\D$ the set of all $k\times k$ matrices $L=(\ell_{ij})$
of nonnegative integers such that the sum of each row and each
column is $n/k$. For any such matrix $L$ observe that there are $
{n!}/({\prod_{i,j} \ell_{ij}!})$ corresponding pairs of balanced
partitions. Therefore,
\begin{eqnarray}\label{eq:secondidentity}
\E Z^2=\sum_{L\in \D} \frac{n!}{\prod_{i=1}^k\prod_{j=1}^k
\ell_{ij}!}\cdot \left[1-\frac{2}{k}+\sum_{i=1}^k\sum_{j=1}^k
\left(\frac{\ell_{ij}}{n}\right)^2\right]^{cn}.
\end{eqnarray}

To get a feel for the sum in~\eqref{eq:secondidentity} observe
that the term corresponding to $\ell_{ij} = n/k^2$ for all $i,j$,
alone, is $\Theta(n^{-(k^2-1)/2}) \cdot [k(1-1/k)^c]^{2n}$. In
fact, the terms corresponding to matrices for which $\ell_{ij} =
n/k^2 \pm O(\sqrt{n})$ already sum to $\Theta((\E Z)^2)$. To
establish $\E Z^2 = O((\E Z)^2)$ we will show that for $c \leq
\h_{k-1}$ the terms in the sum (2) decay exponentially in their
distance from $(\ell_{ij})=(n/k^2)$ and apply
Lemma~\ref{lem:laplacebirkhoff} below. This lemma is a variant of
the classical Laplace method of asymptotic analysis in the case of
the Birkhoff polytope ${\cal B}_k$, \ie the set of all $k \times
k$ doubly stochastic matrices. For a matrix $A\in{\cal B}_k$ we
denote by $\rho_A$ the square of its 2-norm, i.e. $\rho_A\equiv
\sum_{i,j} a_{ij}^2=\|A\|_{2}^2$. Moreover, let $\mathcal H(A)$
denote the entropy of $A$, which is defined as
\begin{eqnarray}\label{eq:def entropy}
{\mathcal{H}}(A)\equiv-\frac{1}{k}\sum_{i=1}^k\sum_{j=1}^ka_{ij}\log
a_{ij}\, .
\end{eqnarray}
Finally, let $J_k \in {\cal B}_k$ be the constant $\frac{1}{k}$
matrix.

\begin{lemma}\label{lem:laplacebirkhoff}
 Assume that $\varphi:{\cal B}_k\to \R$ and $\beta>0$ are such that for every $A\in {\cal B}_k${\rm ,}
$$
{\mathcal{H}}(A)+\varphi(A)\le
{\mathcal{H}}(J_k)+\varphi(J_k)-\beta(\rho_A-1) \enspace .
$$
Then there exists a constant $C = C(\beta,k)>0$ such that
\begin{eqnarray}\label{eq:sum} \sum_{L\in \D} \frac{n!}{\prod_{i=1}^k\prod_{j=1}^k
\ell_{ij}!}\cdot\exp \left[n \cdot
\varphi\left(\frac{k}{n}L\right)\right]\le
%\left[
%\left(\frac{2}{\beta}\right)^{(k-1)^2/2}\cdot\frac{2k^{2k-1}}{(\pi
%n)^{k-1}}e^{\frac{k^4\beta}{4n}}+3n^{k^2}e^{-\frac{n\beta}{4k^2}}\right]
\frac{C}{n^{k-1}} \cdot \left(k^2e^{\varphi(J_k)}\right)^n.
\end{eqnarray}
\end{lemma}
\vglue9pt

 The proof of Lemma~\ref{lem:laplacebirkhoff} is
presented in Section~\ref{section:appendix}.

Let $\St_k$ denote the set of all $k\times k$
\zima{row-stochastic} matrices. For $A \in \St_k$ define
\begin{eqnarray*}
g_c(A)&=&-\frac{1}{k}\sum_{i=1}^k\sum_{j=1}^ka_{ij}\log a_{ij}+c\log
\left(1-\frac{2}{k}+\frac{1}{k^2}\sum_{i=1}^k\sum_{j=1}^ka_{ij}^2\right)\\
&\equiv&
{\mathcal{H}}(A)+c \, {\mathcal{E}}(A).
\end{eqnarray*}

The heart of our analysis is the following inequality. Recall that
$\h_{k-1}=\break(k-1) \log (k-1)$.
\begin{theorem}\label{thm:inequality} For every $A \in \St_k$  and $c \le \h_{k-1}${\rm ,} $
g_c(J_k)\ge g_c(A)$.
\end{theorem}

Theorem~\ref{thm:inequality} is a consequence of a general
optimization principle that we will prove in
Section~\ref{section:simplex} and which is of independent
interest. We conclude this section by showing how
Theorem~\ref{thm:inequality} implies $\ex Z^2 = O((\ex Z)^2)$ and,
thus, Proposition~\ref{thm:threshold}.

For any $A \in {\mathcal B}_k \subset \St_k$ and $c<\h_{k-1}$ we
have
\begin{eqnarray*}
g_c(J_k)-g_c(A)&= &g_{\h_{k-1}}(J_k)-g_{\h_{k-1}}(A)+(\h_{k-1}-c)
\log\left(1+\frac{\rho_A-1}{(k-1)^2}\right) \\
&\geq& (\h_{k-1}-c)
\frac{\rho_A-1}{2(k-1)^2} \enspace ,
\end{eqnarray*}
where for the inequality we applied Theorem~\ref{thm:inequality}
with $c=\h_{k-1}$ and used that $\rho_A \leq k$ so that
$\frac{\rho_A-1}{(k-1)^2}\leq \frac12$. Thus,   for every
$c<\h_{k-1}$ and every $A \in {\mathcal B}_k$
    \begin{equation} \label{expo}
    g_c(A) \leq g_c(J_k)-\frac{\h_{k-1}-c}{2(k-1)^2} \cdot (\rho_A-1) \enspace .
    \end{equation} 
Setting $\beta=(\h_{k-1}-c)/(2(k-1)^2)$ and applying  
Lemma~\ref{lem:laplacebirkhoff} with $\varphi(\cdot)= c\,
{\mathcal{E}}(\cdot)$ yields $\ex Z^2 = O((\ex Z)^2)$.

One can interpret the maximization of $g_c$ geometrically by
recalling that the vertices of the
Birkhoff polytope are the $k!$ permutation matrices (each such matrix
having one non-zero element in each row and column) and $J_k$ is its
barycenter.  By
convexity, $J_k$ is the {\em maximizer\/} of the entropy over
$\mathcal{B}_k$ and the {\em minimizer\/} of the 2-norm. By the
same token, the permutation matrices are {\em minimizers\/} of the
entropy and {\em maximizers\/} of the 2-norm. The constant $c$ is,
thus, the control parameter determining the relative importance of
each quantity. Indeed, it is not hard to see that for sufficiently
small $c$, $g_c$ is maximized by $J_k$ while for sufficiently
large $c$ it is not. The pertinent question is when does the
transition occur, \zima{\ie what is the smallest value of $c$ for
which the norm gain away from $J_k$ makes up for the entropy
loss.} Probabilistically, this is the point where the second
moment explodes (relative to the square of the expectation), as
the dominant contribution stops corresponding to
\mbox{uncorrelated $k$-colorings, \ie to $J_k$.}

The generalization from ${\cal B}_k$ to $\St_k$ is motivated by
the desire to exploit the product structure of  
 the polytope
$S_k$ and Theorem 7 is optimal with respect to $c$, up to an additive
constant. At the same time, it is easy to see that the maximizer of
$g_c$ over $\cal B_k$ is {\em not} $J_k$ already when $c=u_k-1$, e.g.\
$g_c(J_k)<g_c(A)$ for $A=\frac{1}{k-1}J_k+\frac{k-2}{k-1}I$. In other
words, applying the second moment method to balanced $k$-colorings
cannot possibly match the first moment upper bound.

\section{Optimization on products of simplices}\label{section:simplex}

In this section we will prove an inequality which is the main step
in the proof of Theorem~\ref{thm:inequality}. This will be done in
a more general framework since the greater generality, beyond its
intrinsic interest, actually leads to a simplification over the
``brute force" argument.

In what follows we denote by $\Delta_k$ the $k$-dimensional
simplex $\{(x_1,\dots ,x_k)\in [0,1]^k:\ \sum_{i=1}^kx_i=1\}$ and
by $S^{k-1}\subset \mathbb{R}^k$ the unit Euclidean sphere
centered at the origin. Recall that $\St_k$ denotes the set of all
$k\times k$ (row) stochastic matrices. For $1 \le \rho\le k$ we denote
by $\St_k(\rho)$ the set of all $k \times k$ stochastic matrices
with 2-norm $\sqrt{\rho}$, \ie $\St_k(\rho)=\left\{A \in\St_k;\
||A||_2^2=\rho\right\}$.

\begin{definition}\label{def:hoolahoop}
For $\frac{1}{k} \le r \le 1$, let $s^*(r)$ be the unique vector
in $\Delta_k$ of the form $(x,y,\dots ,y)$ having 2-norm
$\sqrt{r}$. Observe that
$$
x=x_r\equiv\frac{1+\sqrt{(k-1)(kr-1)}}{k} \qquad \mbox{and} \qquad
y=y_r\equiv \frac{1-x_r}{k-1} \enspace .
$$
Given $h:[0,1]\to \R$ and an integer $k>1$ we define a function
$f:[1/k,1]\to \R$ as
    \begin{equation} \label{eq:f_def}
    f(r) = h\left(x_r\right)+(k-1) \cdot h\left(y_r\right) \enspace .
    %h\left(\frac{1+\sqrt{(k-1)(kr-1)}}{k}\right)+(k-1) \cdot h\left(\frac{k-1-\sqrt{(k-1)(kr-1)}}{k(k-1)}\right) \enspace .
    \end{equation} 
\end{definition}

Our main inequality provides a sharp bound for the maximum of
entropy-like functions over stochastic matrices with a given
2-norm. In particular, in Section~\ref{sec:asympt} we will prove
Theorem~\ref{thm:inequality} by applying Theorem
\ref{thm:stochastic} below to the function $h(x)=-x\log x$.

\begin{theorem}\label{thm:stochastic} Fix an integer $k> 1$ and let
$h:[0,1]\to \mathbb{R}$ be a continuous strictly concave function{\rm ,}
which is six times differentiable on $(0,1)$. Assume that
$h'(0^+)=\infty${\rm ,} $h'(1^-)>-\infty$ and $h^{(3)}>0${\rm ,} $h^{(4)}<0${\rm ,}
$h^{(6)}<0$ point-wise. Given $1\leq\rho\leq k${\rm ,} for $A\in
\St_k(\rho)$ define
$$
H(A)=\sum_{i=1}^k\sum_{j=1}^kh(a_{ij}).
$$
Then{\rm ,} for $f$ as in~{\rm \eqref{eq:f_def},}
    \begin{equation} \label{eq:bmro}
H(A)\le \max\left\{ m\cdot k \, h\left( \frac{1}{k} \right) +(k-m)
\cdot f \left( \frac{k\rho-m}{k(k-m)} \right);\ 0\le m\le
\frac{k(k-\rho)}{k-1} \right\}.
    \end{equation} 
\end{theorem}

To understand the origin of the right-hand side
in~\eqref{eq:bmro}, consider the following. Given $1 \le \rho \le
k$ and an integer $0\le m\le \frac{k(k-\rho)}{k-1}$, let
$B_{\rho}(m) \in \St_k(\rho)$ be the matrix whose first $m$ rows
are the constant $1/k$ vector and the remaining $k-m$ rows are the
vector $s^*\left(\frac{k\rho-m}{k(k-m)}\right)$. Define
$Q_{\rho}(m) = H(B_{\rho}(m))$. Theorem~\ref{thm:stochastic} then
asserts that $H(A) \leq{\max_{m}} Q_{\rho}(m)$, where $0\le m\le
\frac{k(k-\rho)}{k-1}$ is real.

To prove Theorem~\ref{thm:stochastic} we observe that if $\rho_i$
denotes the squared 2-norm of the $i$-th row then
\begin{equation}\label{eq:indulge}
\max_{A \in \St_k(\rho)} H(A) \;\;= \max_{(\rho_1,\dots ,\rho_k)
\in \rho \Delta_k} \; \sum_{i=1}^k \max \left\{\hat{h}(s); s \in
\Delta_k \cap \sqrt{\rho_i} \, S^{k-1}\right\} \enspace ,
\end{equation}
where $\hat h(s)=\sum_{j=1}^k h(s_j)$. %That is, we can prescribe
%the sum of the squared elements of each row, maximize $H(\cdot)$
%subject to these conditions, and then maximize over
%$(\rho_1,\dots ,\rho_k) \in \rho \Delta_k$.
The crucial point, reflecting the product structure of $\St_k$, is
that to maximize the sum in~\eqref{eq:indulge} it suffices to
maximize $\hat h$ in each row  independently. The maximizer of
each row is characterized by the following proposition:

\begin{proposition}\label{prop:vector} 
Fix an integer $k\ge 1$ and let $h:[0,1]\to \mathbb{R}$ be a continuous
strictly concave function which is three times differentiable on
$(0,1)$. Assume that $h'(0^+)=\infty${\rm ,} 
and $h'''>0$ point-wise. Fix $\frac{1}{k}\le r\le 1$ and assume
that $s=(s_1,\dots ,s_k)\in \Delta_k\cap (\sqrt{r} \, S^{k-1})$ is
such that
$$
\hat h(s)\equiv\sum_{i=1}^k h(s_i)=\max\left\{\sum_{i=1}^k
h(t_i);\ (t_1,\dots ,t_k)\in \Delta_k\cap \sqrt{r} \,
S^{k-1}\right\}.
$$
Then{\rm ,} up to a permutation of the coordinates{\rm ,} $s=s^*(r)$ where
$s^*(r)$ is as in Definition~{\rm \ref{def:hoolahoop}.}
\end{proposition}

Thus, if $\rho_i$ denotes the squared 2-norm of the $i$-th row of
$A\in \St_k$,  Proposition~\ref{prop:vector} implies that $H(A)\le
F(\rho_1,\dots ,\rho_k)\equiv \sum_{i=1}^k f\left(\rho_i\right)$,
where $f$ is as in (\ref{eq:f_def}). Hence, to prove
Theorem~\ref{thm:stochastic} it suffices to give an upper bound on
$F(\rho_1,\dots ,\rho_k)$, where $(\rho_1,\dots ,\rho_k) \in
\rho\Delta_k \cap [1/k,1]^k$. This is another optimization problem
on a symmetric polytope and had $f$ been concave it would be
trivial. Unfortunately, in general, $f$ is not concave (in
particular, it is not concave when $h(x)=-x \log x$).
Nevertheless, the conditions of Theorem~\ref{thm:stochastic} on
$h$ suffice to %deduce
\zima{impart} some properties \zima{on} $f$:
\begin{lemma}\label{lem:f} Let $h:[0,1]\to \R$ be six times
differentiable on $(0,1)$ such that $h^{(3)}>0${\rm ,} $h^{(4)}<0$ and
$h^{(6)}<0$ point-wise. Then the function $f$ defined
in~{\rm \eqref{eq:f_def}} satisfies $f^{(3)}<0$ point-wise.
\end{lemma}

The following lemma is the last ingredient in the proof of
Theorem~\ref{thm:stochastic} as it will allow us to make use of
Lemma~\ref{lem:f} to bound $F$.

\begin{lemma}\label{lem:ineqf}
Let $\psi:[0,1]\to \R$ be continuous on $[0,1]$ and three times
differentiable on $(0,1)$. Assume that $\psi'(1^-)=-\infty$ and
$\psi^{(3)}<0$ point-wise. Fix $\gamma \in (0,k]$ and let
$s=(s_1,\dots ,s_k)\in [0,1]^k\cap\gamma \Delta_k$. Then
$$
\Psi(s) \equiv \sum_{i=1}^k \psi(s_i) \leq
\max\left\{m\psi(0)+(k-m)\psi\left(\frac{\gamma}{k-m}\right);\ m
\in[0,k-\gamma] \right\}.
$$
\end{lemma}

%Indeed, by Proposition~\ref{prop:vector} and the definition of
%$f$, for every $i\in \{1,\dots , k\}$
%\begin{eqnarray*}
%\sum_{j=1}^k h(a_{ij}) & \le &
%h\left(\frac{1+\sqrt{(k-1)(k\rho_i-1)}}{k}\right)+(k-1)h\left(\frac{k-1-\sqrt{(k-1)(k\rho_i-1)}}{k(k-1)}\right)
%\\
%& = & f\left(\rho_i\right) .
%\end{eqnarray*}

To prove Theorem~\ref{thm:stochastic} we define $\psi:[0,1]\to \R$
as $\psi(x)=f\left(\frac{1}{k}+\frac{k-1}{k}x\right)$.
Lemma~\ref{lem:f} and our assumptions on $h$  imply that $\psi$
satisfies the conditions of Lemma \ref{lem:ineqf} (the assumption
that $h'(0^+)=\infty$ implies that $\psi'(1^-)=-\infty$). Hence,
applying Lemma \ref{lem:ineqf} with $\gamma=\frac{k(\rho-1)}{k-1}$
yields Theorem~\ref{thm:stochastic}, \ie
\begin{eqnarray*}
F(A)& = &\sum_{i=1}^k \psi\left(\frac{k\rho_i-1}{k-1}\right) \\
&\le&
\max\left\{m\,
\psi(0)+(k-m)\psi\left(\frac{k(\rho-1)}{(k-1)(k-m)}\right);  \ m\!
\in\!\left[0,k-\frac{k(\rho-1)}{k-1}\right]\! \right\}\! .
\end{eqnarray*}

\Subsec{ Proof of Proposition~{\rm \ref{prop:vector}}}
When $r=1$ there is nothing to prove, so assume that $r<1$. We
begin by observing that $s_i>0$ for every $i\in \{1,\dots , k\}$.
Indeed, for the sake of contradiction, we may assume without loss of
generality (since $r<1$) that $s_1=0$ and $s_2\geq s_3>0$. Fix
$\e>0$ and set
$$
\mu(\e)=\frac{s_2-s_3+\e-\sqrt{(s_2-s_3-\e)^2+4\e(s_3-\e)}}{2}\quad
\mathrm{and}\quad \nu(\e)=-\mu(\e)-\e \enspace.
$$
\zima{Let} $v(\e)=(\e,s_2+\mu(\e),s_3+\nu(\e),s_4,\dots ,s_k)$.
\zima{Our} choice of $\mu(\e)$ and $\nu(\e)$ ensures that for $\e$
small enough $v(\e) \in \Delta_k\cap (\sqrt{r}\cdot S^{k-1})$.
Recall that, by assumption, $h'(0)=\infty$ and $h'(x)<\infty$ for
$x\in(0,1)$. When $s_2 > s_3$ it is clear that $|\mu'(0)| <
\infty$ and, thus, $\left.\frac{d}{d\e}\hat
h(v(\e))\right|_{\e=0}=\infty$. On the other hand, when
$s_2=s_3=s$ it is not hard to see that $$\left.\frac{d}{d\e}\hat
h(v(\e))\right|_{\e=0}=h'(0^{+}) -h'(s)+sh^{''}(s) = \infty.
$$
Thus, in both cases, we have $\left.\frac{d}{d\e}\hat
h(v(\e))\right|_{\e=0}=\infty$ which contradicts the maximality of
$\hat h(s)$.

Since $s_i>0$ for every $i$ \dima{(and, therefore, $s_i<1$ as
well)}, we may use Lagrange multipliers to deduce that there are
$\lambda,\mu \in \dima{\mathbb{R}}$ such that for every $i\in
\{1,\dots , k\}$, $h'(s_i)=\lambda s_i+\mu$. Observe that if we
let $\psi(u)=h'(u)-\lambda u$ then $\psi''=h'''>0$, \ie $\psi$ is
strictly convex. It follows in particular that
$|\psi^{-1}(\mu)|\le 2$. Thus, up to a permutation of the
coordinates, we may assume that there is an integer $1\le m\le k$
and $a,b\in \dima{(0,1)}$ such that $s_i=a$ for $i\in \{1,\dots ,
m\}$ and $s_i=b$ for $i\in \{m+1,\dots , k\}$. Without loss of
generality $a\ge b$ (so that in particular $a\ge 1/k$ and $b\le
1/k$). Since $ma+(k-m)b=1$ and $ma^2+(k-m)b^2=r$, it follows that
$$
a=\frac{1}{k}+\frac{1}{k}\sqrt{\frac{k-m}{m}(kr-1)}\quad
\mathrm{and}\quad
b=\frac{1}{k}-\frac{1}{k}\sqrt{\frac{m}{k-m}(kr-1)} \enspace .
$$
(The choice of the minus sign in the solution of the quadratic
equation defining $b$ is correct since $b\le 1/k$.) Define
$\alpha,\beta:[1,r^{-1}]\to \R$ by
$$
\alpha(t)=\frac{1}{k}+\frac{1}{k}\sqrt{\frac{k-t}{t}(kr-1)} \quad
\mathrm{and}\quad
\beta(t)=\frac{1}{k}-\frac{1}{k}\sqrt{\frac{t}{k-t}(kr-1)}
\enspace .
$$
Furthermore, set $ \varphi(t)=t\cdot h(\alpha(t))+ (k-t)\cdot
h(\beta(t))$, so that $\hat h(s)=\varphi(m)$.

The proof will be complete once we check that $\varphi$ is
strictly decreasing. Observe that
\begin{eqnarray*}
 t\alpha(t)+(k-t)\beta(t) & = & 1 \\
t\alpha(t)^2+(k-t)\beta(t)^2 & = & r\enspace .
\end{eqnarray*}
Differentiating these identities we find that
\begin{eqnarray*}
\alpha(t)+t\alpha'(t)-\beta(t)+(k-t)\beta'(t) & = & 0 \\
\alpha(t)^2+2t\alpha(t)\alpha'(t)-\beta(t)^2+2(k-t)\beta(t)\beta'(t)&
=& 0 \enspace ,
\end{eqnarray*}
implying
$$
\alpha'(t)=-\frac{\alpha(t)-\beta(t)}{2t}\quad \mathrm{and}\quad
\beta'(t)=-\frac{\alpha(t)-\beta(t)}{2(k-t)} \enspace .
$$
Hence,
\begin{eqnarray*}
\varphi'(t)&=&h(\alpha(t))-h(\beta(t))\dima{+t\alpha'(t)h'(\alpha(t))+}(k-t)\beta'(t)h'(\beta(t))\\
&=& h(\alpha(t))-h(\beta(t))-
\frac{\alpha(t)-\beta(t)}{2}[h'(\alpha(t))+h'(\beta(t))] \enspace
.
\end{eqnarray*}

Therefore, in order to show that $\varphi'(t)<0$, it is enough to
prove that if $0\le \beta<\alpha<1$ then
$$
h(\alpha)-h(\beta)-\frac{\alpha-\beta}{2}[h'(\alpha)+h'(\beta)]<0
\enspace .
$$
Fix $\beta$ and define $\zeta:[\beta,1]\to \R$ by
$\zeta(\alpha)=h(\alpha)-h(\beta)-\frac{\alpha-\beta}{2}[h'(\alpha)+h'(\beta)]$.
Now,
$$
\zeta'(\alpha)=\frac{\alpha-\beta}{2}\left(\frac{h'(\alpha)-h'(\beta)}{\alpha-\beta}-h''(\alpha)\right).
$$
By the Mean Value Theorem there is $\beta<\theta<\alpha$ such that
$$
\zeta'(\alpha)=\frac{\alpha-\beta}{2}[h''(\theta)-h''(\alpha)]<0,
$$
since $h'''>0$. This shows that $\zeta$ is strictly decreasing.
Since $\zeta(\beta)=0$ it follows that for $\alpha\in (\beta,1]$,
$\zeta(\alpha)<0$, which concludes the proof of
Proposition~\ref{prop:vector}.\hfill\qed

\Subsec{ Proof of Lemma~{\rm \ref{lem:f}}}
If we make the linear change of variable $z=(k-1)(kx-1)$ then our
goal is to show that the function $g:[0,(k-1)^2]\to \R$, given by
$$
g(z)=h\left(\frac{1}{k}+\frac{\sqrt{z}}{k}\right)+(k-1)h\left(\frac{1}{k}-\frac{\sqrt{z}}{k(k-1)}\right),
$$
satisfies $g'''<0$ point-wise. Differentiation gives
\begin{eqnarray*}
8kz^{5/2}g'''(z)&=&\frac{z}{k^2}
\left[h'''\left(\frac{1}{k}+\frac{\sqrt{z}}{k}\right)-
\frac{1}{(k-1)^2}h'''\left(\frac{1}{k}-\frac{\sqrt{z}}{k(k-1)}\right)\right]\\
& & -\ \frac{3\sqrt{z}}{k}\left[
h''\left(\frac{1}{k}+\frac{\sqrt{z}}{k}\right)+\frac{1}{k-1}h''\left(\frac{1}{k}-\frac{\sqrt{z}}{k(k-1)}\right)\right]\\
& &  +\ 3\left[h'\left(\frac{1}{k}+\frac{\sqrt{z}}{k}\right)-
h'\left(\frac{1}{k}-\frac{\sqrt{z}}{k(k-1)}\right)\right].
\end{eqnarray*}
Denote $a=\frac{\sqrt{z}}{k}$ and $b=\frac{\sqrt{z}}{k(k-1)}$.
Then $8kz^{5/2}g'''(z)=\psi(a)-\psi(-b)$, where
$$
\psi(t)=t^2h'''\left(\frac{1}{k}+t\right)-3th''\left(\frac{1}{k}+t\right)+3h'\left(\frac{1}{k}+t\right).
$$
Now
$$
\psi'(t)=t^2h''''\left(\frac{1}{k}+t\right)-th'''\left(\frac{1}{k}+t\right).
$$
The assumptions on \dima{$h'''$ and $h''''$} imply that
$\psi'(t)<0$ for $t>0$, and since $a\ge b$, it follows that
$\psi(a)\le \psi(b)$. Since
$8kz^{5/2}g'''(z)=\psi(a)-\psi(-b)=\big[\psi(a)-\psi(b)\big]+\big[\psi(b)-\psi(-b)\big]$,
it suffices to show that for every $b>0$,
$\zeta(b)=\psi(b)-\psi(-b)<0$. Since $\zeta(0)=0$, this will
follow once we verify that $\zeta'(b)<0$ for $b>0$. Observe now
that $\zeta'(\beta) = b \chi(b)$, where
$$
\chi(b)=b\left[h''''\left(\frac{1}{k}+b\right)+h''''\left(\frac{1}{k}-b\right)\right]-
\left[h'''\left(\frac{1}{k}+b\right)-h'''\left(\frac{1}{k}-b\right)\right].
$$
Our goal is to show that $\chi(b)<0$ for $b>0$, and since
$\chi(0)=0$ it is enough to show that $\chi'(b)<0$. But
$$
\chi'(b)=b\left[h^{(5)}\left(\frac{1}{k}+b\right)-h^{(5)}\left(\frac{1}{k}-b\right)\right],
$$
so that the required result follows from the fact that $h^{(5)}$
is strictly decreasing.

\Subsec{ Proof of Lemma~{\rm \ref{lem:ineqf}}}
Before proving Lemma~\ref{lem:ineqf} we require one more
preparatory fact.

\begin{lemma}\label{lem:abm} Fix $0<\gamma<k$. Let $\psi:[0,1]\to \R$ be continuous on $[0,1]$ and three times
differentiable on $(0,1)$. Assume that $\psi'(1^-)=-\infty$ and
$\psi'''<0$ point-wise. Consider the set $A\subset \R^3$ defined
by
$$
A=\left\{(a,b,\ell)\in \dima{(}0,1]\times [0,\dima{1]} \times
\dima{(0},k];\ b<a \ \mathrm{and}\ \ell
a+(k-\ell)b=\gamma\right\}.
$$
Define $g:A\to \R$ by $g(a,b,\ell)=\ell \psi(a)+(k-\ell)\psi(b)$.
\dima{If $(a,b,\ell)\in A$ is such that
$g(a,b,\ell)=\max_{(a,b,\ell)\in A}g(a,b,\ell)$ then $a =
\gamma/\ell$}.
\end{lemma}

{\it Proof of  Lemma}~\ref{lem:abm}.
Observe that if \dima{$b=0$ or $\ell = k$} we are done. Therefore,
assume that \dima{$b>0$ and $\ell<k$}. We claim that $a<1$.
Indeed, if $a=1$ then $b=\frac{\gamma-\ell}{k-\ell}<1$, implying
that for small enough $\e>0$,  $w(\e)\equiv\left(1-\e,
b+\frac{\ell\e}{k-\ell},\ell\right)\in A$. But
$\left.\frac{d}{d\e}g(w(\e))\right|_{\e=0}=-\ell \psi'(1^-)+\ell
\psi'(b)=\infty$, which contradicts the maximality of
$g(a,b,\ell)$.

\dima{Since $a \in (0,1)$ and $\ell \in (0,k)$} we can use
Lagrange multipliers to deduce that there is $\lambda\in
\dima{\R}$ such that $\ell \psi'(a)=\lambda \ell$,
$(k-\ell)\psi'(b)=\lambda(k-\ell)$ and
$\psi(a)-\psi(b)=\lambda(a-b)$. Combined, these imply
$$
\psi'(a)=\psi'(b)=\frac{\psi(a)-\psi(b)}{a-b} \enspace . \pagebreak
$$
\dima{By} the Mean Value Theorem, there exists $\theta\in (b,a)$
such that  \dima{$\psi'(\theta)=\frac{\psi(a)-\psi(b)}{a-b}$. But,
since $\psi'''<0$, $\psi'$ cannot take the same value three times,
yielding the desired contradiction.}
\Endproof\vskip4pt

We now turn to the proof of Lemma~\ref{lem:ineqf}. Let $s \in
[0,1]^k\cap \gamma\Delta_k$ be such that $\Psi(s)$ is maximal. If
$s_1=\cdots=s_k=1$ then we are done, so we assume that there
exists $i$ for which $s_i<1$. Observe that in this case $s_i<1$
for every $i\in \{1,\dots , k\}$. Indeed, assuming the contrary we
may also assume without loss of generality that $s_1=1$ and
$s_2<1$. For every $\e>0$ consider the vector
$u(\e)=(1-\e,s_2+\e,s_3,\dots ,s_k)$. For $\e$ small enough
$u(\e)\in [0,1]^k\cap\gamma\Delta_k$. But
$\left.\frac{d}{d\e}\Psi(u(\e))\right|_{\e=0}=\infty$, which
contradicts the maximality of $\Psi(s)$.

Without loss of generality we can further assume that
$s_1,\dots ,s_q>0$ for some \dima{$q\leq k$} and $s_i=0$ for all
$i>q$. Consider the function $\tilde{\Psi}(t) =
\sum_{i=1}^q\psi(t_i)$ defined on $[0,1]^q\cap \gamma\Delta_q$.
Clearly, $\tilde{\Psi}$ is maximal at $(s_1,\dots ,s_q)$. Since
\dima{$s_i\in(0,1)$} for every $i\in \{1,\dots , q\}$, we may use
Lagrange multipliers to deduce that there is $\lambda\in
\dima{\R}$ such that for every $i\in \{1,\dots , q\}$,
$\psi'(s_i)=\lambda$. Since $\psi'''<0$, $\psi'$ is strictly
concave. It follows in particular that the equation
$\psi'(y)=\lambda$ has at most two solutions, so that up to a
permutation of the coordinates we may assume that there is an
integer $0\le \ell\le q$ and \dima{$0\leq b<a\leq 1$} such that
$s_i=a$ for $i\in \{1,\dots , \ell\}$ and $s_i=b$ for $i\in
\{\ell+1,\dots , q\}$. Now, using the notation of Lemma
\ref{lem:abm} we have that $(a,b,\ell)\in A$ so that
\begin{eqnarray*}
\Psi(s)&=&(k-q)\psi(0) + g(a,b,\ell)\\[5pt]
&\le& (k-q)\psi(0) + \max\left\{\theta
\psi(0)+(q-\theta)\psi\left(\frac{\gamma}{q-\theta}\right);\
\theta \in[0,q-\gamma] \right\} \\[5pt]
& \leq &
\zima{\max\left\{m\psi(0)+(k-m)\psi\left(\frac{\gamma}{k-m}\right);\
m \in[0,k-\gamma] \right\}} \enspace .
\end{eqnarray*}
%This is precisely the assertion of Lemma~\ref{lem:ineqf} (letting
%$m=k-q+\theta$).

\section{Proof of Theorem~\ref{thm:inequality}}\label{sec:asympt}

Let $h(x)=-x\log x$ and note that $h'(x)=-\log x-1$,
$h'''(x)=\frac{1}{x^2}$, $h^{(4)}(x)=\frac{-2}{x^3}$ and
$h^{(6)}(x)=\frac{-24}{x^5}$, so that the conditions of Theorem
\ref{thm:stochastic} are satisfied in this particular case. By
Theorem \ref{thm:stochastic} it is, thus, enough to show that for
$c \leq \h_{k-1} = (k-1) \log (k-1)$,
\begin{multline}\label{eq:goal}
\frac{m\log
k}{k}+\frac{k-m}{k}f\left(\frac{k\rho-m}{k(k-m)}\right)+c\log\left(1-\frac{2}{k}
+\frac{\rho}{k^2}\right)\\[5pt]
\le
 \log k+2c\log\left(1-\frac{1}{k}\right),
\end{multline}
for every \pagebreak $1\le \rho\le k$ and $0\le m\le \frac{k(k-\rho)}{k-1}$.
Here $f$ is as in (\ref{eq:f_def}) for $h(x)=-x\log x$. Inequality
\eqref{eq:goal} simplifies to
\begin{eqnarray}\label{eq:with m}
c\log\left(1+\frac{\rho-1}{(k-1)^2}\right)\le
\left(1-\frac{m}{k}\right)\left[ \log
k-f\left(\frac{k\rho-m}{k(k-m)}\right)\right].
\end{eqnarray}
Setting $t=m/k$, $s=\rho-1$ and using the inequality $\log(1+a)\le
a$, it suffices to demand that for every $0\le t\le
1-\frac{s}{k-1}$ and $0\le s\le k-1$,
\begin{eqnarray}\label{eq:goodgoal}
 \frac{cs}{(k-1)^2}\le (1-t)\left[
f\left(\frac{1}{k}\right)-f\left(\frac{1}{k}+\frac{s}{k(1-t)}\right)\right].
\end{eqnarray}

To prove \eqref{eq:goodgoal} we define $\eta : (0,1-1/k] \to
\reals$ by
$$
\eta(y)=\frac{f\left(\frac{1}{k}\right)-f\left(\frac{1}{k}+y\right)}{y}
\enspace ,
$$
and $\eta(0)=-f'\left(\frac{1}{k}\right)=\frac{k}{2}$, making
$\eta$ continuous on $[0,1-1/k]$. Observe that (\ref{eq:goodgoal})
reduces to
\begin{eqnarray*}\label{eq:better}
c\le \frac{(k-1)^2}{k}\cdot \eta\left(\frac{s}{k(1-t)}\right).
\end{eqnarray*}
Now, $\eta'(y)=\frac{\zeta(y)}{y^2}$, where
$\zeta(y)=f\left(\frac{1}{k}+y\right)-f\left(\frac{1}{k}\right)-yf'\left(\frac{1}{k}+y\right)$.
Observe that $\zeta'(y)=-yf''\left(\frac{1}{k}+y\right)$ so, by
Lemma \ref{lem:f}, $\zeta$ can have at most one zero in
$\left(0,1-\frac{1}{k}\right)$. A straightforward computation
gives that $\zeta\left(\frac{(k-2)^2}{k(k-1)}\right)=0$, so $\eta$
achieves its global minimum on $ \left[0,1-\frac{1}{k}\right]$ at %the points
$y \in \left\{0,\frac{(k-2)^2}{k(k-1)},1-\frac{1}{k}\right\}$.
Direct computation gives
$\eta\left(1-\frac{1}{k}\right)=\frac{k}{k-1}\cdot \log k$,
$\eta\left(\frac{(k-2)^2}{k(k-1)}\right)= \frac{k-1}{k-2}\cdot
\log(k-1)$ and, by definition, $ \eta(0) =\frac{k}{2}$. Hence
\begin{eqnarray}&&\label{eq:never_use}\\
\frac{(k-1)^2}{k}\cdot \eta\left(\frac{s}{k(1-t)}\right)&\ge&
\frac{(k-1)^2}{k}\cdot\min\left\{\frac{k}{2},\frac{k-1}{k-2}\cdot
\log(k-1),\frac{k}{k-1}\cdot \log k\right\}\nonumber
\\
& = &\frac{(k-1)^3}{k(k-2)}\cdot \log(k-1)  > \h_{k-1} \enspace ,
\nonumber
\end{eqnarray}
where~\eqref{eq:never_use} follows from elementary calculus.

\demo{Remark} The above analysis shows that
Theorem~\ref{thm:inequality} is asymptotically optimal. Indeed,
let $A$ be the stochastic matrix whose first $k-1$ rows are the
constant $1/k$ vector and whose last row is the vector $s^*(r)$,
defined in Definition~\ref{def:hoolahoop}, for
$r=\frac{1}{k}+\frac{(k-2)^2}{k(k-1)}$. This matrix corresponds to
$m=k-1$ and $\rho=1+\frac{(k-2)^2}{k(k-1)}$ in~\eqref{eq:with m},
and a direct computation shows that any $c$ for which\break\vglue-11pt\noindent
Theorem~\ref{thm:inequality} holds must satisfy $c < \h_{k-1}+1$.

\section{Appendix: Proof of
Lemma~\ref{lem:laplacebirkhoff}}\label{section:appendix}

If $(\ell_{ij})$ are nonnegative integers such that
$\sum_{i,j}\ell_{ij}=n$, standard Stirling approximations imply
\begin{eqnarray}\label{eq:nozero}
 \frac{n!}{\prod_{i=1}^k\prod_{j=1}^k
\ell_{ij}!}&\le&
\left[\prod_{i=1}^k\prod_{j=1}^k\left(\frac{\ell_{ij}}{n}\right)^{-\ell_{ij}/n}\right]^n
\\
&&\cdot \min \left\{3\sqrt{n},\left[(2\pi
n)^{k^2-1}\prod_{i=1}^k\prod_{j=1}^k
\frac{\ell_{ij}}{n}\right]^{-1/2}\right\}  \enspace .\nonumber
\end{eqnarray}

Since $|\D|\le (n+1)^{(k-1)^2}$, the contribution to the sum in
(\ref{eq:sum}) of the terms for which $\rho_{\frac{k}{n}L} >
1+1/(4k^2)$ can, thus, be bounded by
\begin{eqnarray}\label{lalakis}
3\sqrt{n}(n+1)^{(k-1)^2}\left(e^{{\mathcal{H}}\left(\frac{k}{n}L\right)+\log
k+\varphi\left(\frac{k}{n}L\right)}\right)^n&\le &3n^{k^2}
\left(k^2e^{\varphi(J_k)}\right)^n \cdot e^{-\frac{\beta n}{4k^2}}
\\&=& O\left(n^{-k^2}\right) \left(k^2e^{\varphi(J_k)}\right)^n . \nonumber
\end{eqnarray}

Furthermore, if $L\in \D$ is such that $\rho_{\frac{k}{n}L}\le
1+\frac{1}{4k^2}$, then for every $1\le i,j\le k$ we have
$$
\left(\frac{k}{n}\ell_{ij}-\frac{1}{k}\right)^2\le
\sum_{s=1}^k\sum_{t=1}^k\left(\frac{k}{n}\ell_{st}-\frac{1}{k}\right)^2=\rho_{\frac{k}{n}L}-1\le
\frac{1}{4k^2} \enspace .
$$
Therefore, for such $L$ we must have $\ell_{ij}\ge n/(2k^2)$ for
every $i,j$. \zima{Therefore, by}~\eqref{eq:nozero},
\eqref{lalakis} we get
\begin{multline}\label{eq:twosum}
\sum_{L\in \D} \frac{n!}{\prod_{i=1}^k\prod_{j=1}^k
\ell_{ij}!}\cdot \exp
\left[n\varphi\left(\frac{k}{n}L\right)\right] \\
\le
 \frac{C(\beta,k)}{ n^{(k^2-1)/2}}\cdot
\left(k^2e^{\varphi(J_k)}\right)^n\cdot \sum_{L\in \D}e^{-\beta
n\left(\frac{k^2}{n^2}\rho_L-1\right)} \enspace .
\end{multline}

Denote by $M_k(\R)$ the space of all $k\times k$ matrices over
$\R$ and let $F$ be the subspace of $M_k(\R)$ consisting of all
matrices $X=(x_{ij})$ for which the sum of each row and each
column is $0$. The dimension of $F$ is $(k-1)^2$. Denote by
$B_\infty$ the unit cube of $M_k(\R)$, i.e. the set of all
$k\times k$ matrices $A=(a_{ij})$ such that $a_{ij}\in [-1/2,1/2]$
for all $1\le i,j\le k$. For $L \in \D$ we define
$T(L)=L-\frac{n}{k}J_k+(F\cap B_\infty)$, \ie the tile $F\cap
B_\infty$ shifted by $L-\frac{n}{k}J_k$.

\begin{lemma}\label{sublem:pajor} For every $L\in \D${\rm ,}
$$
e^{-\beta n \left(\frac{k^2}{n^2}\rho_L-1\right)}\le
e^{\frac{k^4\beta}{4n}}\cdot \int_{T(L)}
e^{-\frac{k^2\beta}{2n}\|X\|_2^2}dX \enspace.
$$
\end{lemma}

\Proof 
By the triangle inequality, we see that for any matrix $X$
\begin{equation}\label{malakies}
\left\|L-\frac{n}{k}J_k\right\|_2^2 \; \geq \;
\frac{1}{2}\|X\|_2^2-\left\|L-\frac{n}{k}J_k-X\right\|_2^2 \; \ge
\;  \frac12\|X\|_2^2-k^2\left\|L-\frac{n}{k}J_k-X\right\|_\infty^2
  .
\end{equation}
Thus, for $X \in T(L)$ we have $\|X\|_2^2\leq
2\left(\frac{k^2}{n^2}\rho_L-1\right)\left(\frac{n}{k}\right)^2+\frac{k^2}{2}$,
since
$\left\|L-\frac{n}{k}J_k\right\|_2^2=\left(\frac{k^2}{n^2}\rho_L-1\right)\left(\frac{n}{k}\right)^2$
and $\left\|L-\frac{n}{k}J_k-X\right\|_\infty^2 \leq \frac14$.
Therefore,
\begin{eqnarray*}
\int_{T(L)}e^{-\frac{k^2\beta}{2n}\|X\|_2^2}dX &\ge& \int_{T(L)}
e^{-\frac{k^4\beta}{4n}}\cdot e^{-\beta
n\left(\frac{k^2}{n^2}\rho_L-1\right)}dX\\
&=&
e^{-\frac{k^4\beta}{4n}}\cdot e^{-\beta
n\left(\frac{k^2}{n^2}\rho_L-1\right)}\mathrm{vol}\left(F\cap
B_\infty\right) \enspace .
\end{eqnarray*}
\zima{It is a theorem of Vaaler~\cite{vaaler} that for any
subspace $E$, $\mathrm{vol}\left(E\cap B_\infty\right) \ge 1$ ,
concluding the proof.}
\Endproof\vskip4pt  

Thus, to bound the second sum in~\eqref{eq:twosum} we apply
Lemma~\ref{sublem:pajor} to get
\begin{eqnarray*}
\sum_{L\in \D }e^{-\beta n \left(\frac{k^2}{n^2}\rho_L-1\right)}
&\le&e^{\frac{k^4\beta}{4n}} \sum_{L\in
\D}\int_{T(L)} e^{-\frac{k^2\beta}{2n}\|X\|_2^2} \, dX\\
&\le& e^{\frac{k^4\beta}{4n}} \int_F
e^{-\frac{k^2\beta}{2n}\|X\|_2^2}
\, dX\\
&=& e^{\frac{k^4\beta}{4n}} \int_{\R^{(k-1)^2}}
e^{-\frac{k^2\beta}{2n}\|X\|_2^2} \, dX\\ &=&
e^{\frac{k^4\beta}{4n}} \left(\frac{2\pi n}{\beta
k^2}\right)^{(k-1)^2/2}
\end{eqnarray*}
where we have used the fact that the interiors of the ``tiles"
$\left\{T(L)\right\}_{L\in\D}$ are disjoint, that the Gaussian
measure is rotationally invariant and that $F$ is $(k-1)^2$
dimensional.
 
\demo{Acknowledgements} We are grateful to Cris Moore for several
inspiring conversations in the early stages of this work.

\def\cprime{$'$} \def\cprime{$'$} \def\cprime{$'$}

\end{document}